\documentclass[12pt]{article}
\usepackage{latexsym,amsmath,amssymb,amsfonts,epsfig,graphicx,cite,psfrag}
\usepackage{eepic,xcolor,colordvi,amscd}
\usepackage{setspace}
\usepackage{algorithm}
\usepackage{algorithmic}
\usepackage{mathrsfs}
\usepackage{amsthm}
\usepackage{enumerate}
\usepackage{enumitem}

\newtheorem{theo}{Theorem}[section]
\newtheorem{ques}[theo]{Question}
\newtheorem{lem}[theo]{Lemma}
\newtheorem{coro}[theo]{Corollary}
\newtheorem{observation}[theo]{Observation}

\newtheorem{conj}[theo]{Conjecture}

\newtheorem{prob}[theo]{Problem}

\newtheorem{claim}[theo]{Claim}

\newcommand\MAD{\hbox{\sc mad}}
\newcommand\ceilgapv{\lceil gap(v)\rceil}

\title{Proper orientations and proper chromatic number}

\author{Yaobin Chen\thanks{Shanghai Center for Mathematical Sciences, Fudan University, Shanghai, 200438, China,
{\tt ybchen21@m.fudan.edu.cn}}
\and
Bojan Mohar\thanks{Department of Mathematics,
Simon Fraser University, Burnaby, BC V5A 1S6, Canada, {\tt mohar@sfu.ca}}~\thanks{Supported in part by the NSERC Discovery Grant R611450 (Canada),
and by the Research Project J1-2452 of ARRS (Slovenia).}%
~~\thanks{On leave from IMFM, Department of Mathematics, University of Ljubljana.}\\
\and
Hehui Wu\thanks{Shanghai Center for Mathematical Sciences, Fudan University,
Shanghai, 200438, China,
{\tt hhwu@fudan.edu.cn}}
~\thanks{Supported in part by National Natural Science Foundation of China grant 11931006, National Key Research and Development Program of China (Grant No. 2020YFA0713200), and the Shanghai Dawn Scholar Program grant 19SG01.}
~
}

\date{\today}

\begin{document}

\maketitle

\begin{abstract}
    The proper orientation number $\Vec{\chi}(G)$ of a graph $G$ is the minimum $k$ such that there exists an orientation of the edges of $G$ with all vertex-outdegrees at most $k$ and such that for any adjacent vertices, the outdegrees are different. Two major conjectures about the proper orientation number are resolved. First it is shown, that $\Vec{\chi}(G)$ of any planar graph $G$ is  at most 14. Secondly, it is shown that for every graph, $\Vec{\chi}(G)$ is at most $O(\frac{r\log r}{\log\log r})+\tfrac{1}{2}\MAD(G)$, where $r=\chi(G)$ is the usual chromatic number of the graph, and $\MAD(G)$ is the maximum average degree taken over all subgraphs of $G$. Several other related results are derived. Our proofs are based on a novel notion of fractional orientations.
\end{abstract}

\section{Introduction}

Borowiecki, Grytczuk and Pil\'{s}niak \cite{BoGrPi12} discovered a beautiful fact that every graph admits an orientation of its edges such that the outdegrees of any two adjacent vertices are different. Such orientations can be interpreted as graph colorings and are now known as \emph{proper orientations}. With this interpretation in mind we define the \emph{proper chromatic number}, also called the \emph{proper orientation number}, $\Vec{\chi}(G)$ of a graph $G$ as the minimum value, taken over all proper orientations of $G$, of the maximum outdegree, $\max\{d^+(v)\mid v\in V(G)\}$. Let us observe that several other papers about proper orientations use the indegree form for this parameter, but the two versions are clearly equivalent.

The original interest in proper orientations came from their connection to the 1-2-3-Conjecture of Karo\'{n}ski, {\L}uczak and Thomason \cite{KaLuTh04}. The first systematic study of the proper orientation number can be found in Ahadi and Dehghan \cite{AD13}. Improved results for bipartite graphs came in a paper by Araujo, Cohen, de Rezende, Havet, and Moura \cite{ACdRHM15}, who proved that every bipartite graph satisfies the following:
\begin{equation}
 \Vec{\chi}(G) \le \left\lfloor \tfrac{1}{2}\Bigl(\Delta(G)+\sqrt{\Delta(G)}\,\Bigr) \right\rfloor + 1.
 \label{eq:bipartite bound}
\end{equation}
Prior to our results, this was still the best general upper bound on $\Vec{\chi}(G)$ of bipartite graphs. In this paper, we provide an improvement that is essentially best possible general upper bound.

Before continuing, we should mention that computing $\Vec{\chi}(G)$ is NP-hard. Moreover, given a graph $G$ and a positive integer $k$, deciding whether $\Vec{\chi}(G)\le k$ is NP-complete, even if $G$ is the line graph of a regular graph \cite{AD13}; or if $G$ is a planar subcubic graph \cite{ACdRHM15}; or if $G$ is planar and bipartite with maximum degree 5 \cite{ACdRHM15}. 
These results show that any insight into the proper orientation number is of interest.

In the case of sparse graph families, one can say a bit more. Papers \cite{KHLMM17,ACdRHM15} and \cite{AHLS16_cacti,ai2020proper} treated trees and outerplanar graphs. The main outcome is that the proper orientation number is bounded on these classes if they satisfy some other conditions (e.g. being triangle-free, or bipartite, or sufficiently connected). The paper by Knox et al.~\cite{KHLMM17} shows that every 3-connected bipartite planar graph $G$ satisfies $\Vec{\chi}(G) \le 5$. This was later improved by Noguchi \cite{No20}, who proved that every bipartite planar graph $G$ with minimum degree 3 satisfies $\Vec{\chi}(G) \le 3$.

However, none of these works was able to solve the most intriguing question from \cite{AHLS16_cacti} whether the proper orientation number of all planar graphs is bounded. One of our main results, see Theorem \ref{th:main1} below, resolves this question in the affirmative. 

The key to our results is another question from \cite{ACdRHM15}. First of all, in view of (\ref{eq:bipartite bound}), Araujo et al.~\cite{ACdRHM15} asked the following.

\begin{prob}\label{prob:bipartite bound}
Does there exists a constant $C$ such that every bipartite graph $G$ satisfies $\Vec{\chi}(G)\le \frac{\Delta(G)}{2}+C$?
\end{prob}

We will answer Problem \ref{prob:bipartite bound} in the affirmative.
In fact, the upper bound depending on $\Delta(G)$ will be improved by using the related quantity called the \emph{Maximum Average Degree} of the graph, $\MAD(G)$, which is defined as the largest average degree of all subgraphs of $G$:
$$
   \MAD(G) = \max_{H\subseteq G} \frac{2|E(H)|}{|V(H)|}.
$$
Of course, the maximum in this definition can be taken over all induced subgraphs only.
This is the basic parameter in \emph{sparsity theory} (see \cite{NeOdM12}) and is a well-known upper bound on the chromatic number of any graph: $\chi(G)\le \lfloor \MAD(G)\rfloor + 1$. On the other hand, $\tfrac{1}{2}\MAD(G)$ gives a general lower bound on $\Vec{\chi}(G)$, since a graph with $\MAD(G)=d$ has a vertex with outdegree at least $\lceil d/2\rceil$ under any orientation of its edges.

Our first main result resolves Problem \ref{prob:bipartite bound}.

\begin{theo}\label{th:main3}
Let $G$ be a bipartite graph. Then 
$$
  \bigl\lceil\tfrac{1}{2}\MAD(G)\bigr\rceil \le \Vec{\chi}(G)\le \bigl\lceil\tfrac{1}{2}\MAD(G)\bigr\rceil + 3.
$$ 
This bound is tight since for each $k$, there exist a bipartite graph $G$ with $\bigl\lceil\tfrac{1}{2}\MAD(G)\bigr\rceil=k$ and $\Vec{\chi}(G) = k+3$.
\end{theo}

Theorem \ref{th:main3} implies various previous results on bipartite planar graphs, for example a result of Araujo et al.~\cite{ACdRHM15}, who proved that every tree has proper orientation number at most 4 (which is best possible). Meanwhile, we obtain a bound for arbitrary bipartite planar graphs.

\begin{coro}\label{cor:main3}
Let $G$ be a bipartite planar graph. Then $\Vec{\chi}(G) \le 5$.
\end{coro}

There is a folklore conjecture that the proper orientation number is bounded on the class of all planar graphs, but no finite upper bound was ever established. By extending the method used in the proof of Theorem \ref{th:main3}, we confirm the conjecture by obtaining the following upper bound for planar graphs.

\begin{theo}\label{th:main1}
Let $G$ be a planar graph. Then $\Vec{\chi}(G) \le 14$. 
\end{theo}

In the proof of the upper bound we use the fact that planar graphs are 4-colorable (the Four-Color Theorem \cite{AH89,RSST97}). Let us remark that even without this result we are able to obtain a (slightly weaker) upper bound by using the simpler fact that planar graphs are 5-colorable. 

The bound of Theorem \ref{th:main1} may not be optimal. Planar graphs with $\Vec{\chi}(G) = 10$ have been constructed by Araujo, Havet, Linhares Sales, and Silva \cite{AHLS16_cacti}, but no example with $\Vec{\chi}(G) > 10$ are known.

For 3-colorable planar graphs we have a stronger bound. 

\begin{theo}\label{th:planar3colorable}
Let $G$ be a $3$-colorable planar graph. Then $\Vec{\chi}(G) \le 11$. Moreover, if $G$ is outerplanar, then $\Vec{\chi}(G) \le 10$.
\end{theo}

Theorem \ref{th:planar3colorable} in particular resolves a problem of Araujo et al.~\cite{AHLS16_cacti}, who conjectured that all outerplanar graphs have bounded $\Vec{\chi}(G)$.
Araujo et al.~\cite{AHLS16_cacti,araujo2021properchordal} and Ai et al.~\cite{ai2020proper} proved that various classes of outerplanar graphs (cactus graphs, triangle-free 2-edge-connected outerplanar graphs, and maximal outerplanar graphs whose inner dual is a path) have bounded $\Vec{\chi}(G)$. As for a lower bound, Araujo et al.~\cite{AHLS16_cacti} found outerplanar graphs having $\Vec{\chi}(G)\ge7$.

The proof for outerplanar graphs uses the fact that these graphs are 3-colorable and their maximum average degree is less than 4. The same proof gives a bound for the proper orientation number of more general graphs. In particular it applies to the superclass of all series-parallel graphs. As the series-parallel graphs are precisely the graphs whose tree-width is at most 2, we have the following result.

\begin{coro}
  Let $G$ be a graph whose tree-width is at most $2$. Then $\Vec{\chi}(G)\le 10$.
\end{coro}

Building on their results about bipartite graphs, Araujo et al.~\cite{ACdRHM15} asked a more general question.

\begin{prob}[Araujo et al.~\cite{ACdRHM15}]
\label{prob:MAD upper bound}
Can $\Vec{\chi}(G)$ be bounded above by a function of $\MAD(G)$?
\end{prob}

Our second main result provides a strong answer to Problem \ref{prob:MAD upper bound}, and resolves another basic question. 

\begin{theo}\label{th:main2}
Let $G$ be a graph whose chromatic number is $r$. Then 
$$\Vec{\chi}(G) = O\Bigl(\frac{r\log r}{\log\log r}\Bigr) + \tfrac{1}{2}\MAD(G).$$
\end{theo}

The constant involved in the $O$-notation in Theorem \ref{th:main2} is small. The precise dependence is given in Section \ref{sect:r-partite} as Theorem \ref{thm:explicit upper bound}.

In particular, since $r=\chi(G)\le \MAD(G)+1$, Theorem \ref{th:main2} provides a resolution of Problem \ref{prob:MAD upper bound}.

\begin{coro}\label{cor:main2}
For a graph $G$ with $\MAD(G)=d$, $\Vec{\chi}(G) = O\bigl(\frac{d\log d}{\log\log d}\bigr)$.
\end{coro}

The paper is organized as follows. In Sections \ref{sect:2} and \ref{sect:3} we give our main tools. Then we apply them to prove the main results, first Theorem \ref{th:main3} on bipartite graphs, then Theorem \ref{th:main2}, and finally Theorems \ref{th:main1} and \ref{th:planar3colorable} about the proper orientation number of planar graphs.


\section{Partial and fractional orientations}
\label{sect:2}

An orientation of a graph is a \emph{$k$-orientation} if every vertex has outdegree at most $k$. It is easy to see that a graph $G$ with a $k$-orientation has $\MAD(G)\le 2k$. Hakimi \cite{hakimi1965degrees} proved the converse statement.

\begin{lem}[Hakimi \cite{hakimi1965degrees}]\label{lem:hakimi}
A graph $G$ admits a $k$-orientation if and only if $\MAD(G)\le 2k$.
\end{lem}

If $G$ is a graph, we can describe any of its orientations by specifying, for each edge $uv\in E(G)$, two values $p(u,v)$ and $p(v,u)$, one of which is 1 and the other one is 0. If $p(u,v)=1$, then we say that the edge $uv$ is \emph{oriented from $u$ to $v$}. If we orient only a subset of the edges, we can have $p(u,v) = p(v,u) = 1$ for those edges $uv$ that are left \emph{unoriented}. Such a function $p$ will be called a \emph{partial orientation}. In this paper we shall use a generalized version of partial orientations where we will allow $p(u,v)$ having any value in $[0,1]$. In fact, we will have the following three possibilities for each edge $uv$:
\begin{itemize}[leftmargin=2cm]
    \item[(PFO1)] An \emph{unoriented edge} $uv$ satisfies $p(u,v) = p(v,u) = 1$.
    \item[(PFO2)] An \emph{oriented edge} $uv$ satisfies $p(u,v) = 1$ and $p(v,u) = 0$. In this case we consider the edge $uv$ as being oriented in the direction from $u$ to $v$.
    \item[(PFO3)] A \emph{fractionally oriented edge} $uv$ satisfies $p(u,v) = \alpha$ and $p(v,u) = 1 - \alpha$ for some $\alpha \in (0,1)$. 
\end{itemize}
If all edges satisfy (PFO1) and (PFO2), then we say that $p$ is a \emph{partial orientation} of $G$. If all edges satisfy (PFO1)--(PFO3), then we say that $p$ is a \emph{partial fractional orientation} of $G$, and if (PFO1) never occurs (all edges are oriented), it is a \emph{(fractional) orientation}.

For a partial fractional orientation $p$ we define the following values for each vertex $v\in V(G)$. The \emph{potential outdegree} $d_p(v)$ of $v$ (with respect to the PFO $p$) is defined as follows:
$$
    d_p(v) = \sum_{u\in N(v)} p(v,u),
$$
where $N(v)$ denotes the set of all neighbors of $v$ in $G$. 
The potential outdegree $d_p(v)$ represents the largest possible outdegree of the vertex $v$, obtained if all unoriented edges incident with $v$ would be oriented out of $v$ (and the fractionally oriented edges are unchanged).
The \emph{outdegree} $d_p^+(v)$ of $v$ is the fractional weight of oriented edges incident with $v$, not counting the unoriented edges. It can be expressed as follows:
$$
    d_p^+(v) = \sum_{u\in N(v)} (1-p(u,v)).
$$

Fractional orientations can be turned into usual orientations such that the outdegree of each vertex is just the ``rounding'' of the fractional outdegree (see \cite{LRS11}). 
This result is stated next just for reference although we are not using it. Instead, we are using our own rounding method that is described after the lemma. 

\begin{lem}\label{lem:rounding}
Let $G$ be a graph and let $p$ be a fractional orientation of the edges of $G$. Then there is an orientation $q$ of the edges of $G$ such that
$$
   d_q^+(u) = \Bigl\lfloor \sum_{v\in N(u)} p(u,v) \Bigr\rfloor \text{ or }
             \Bigl\lceil \sum_{v\in N(u)} p(u,v) \Bigr\rceil.
$$

\end{lem}

Recall that if $p$ is a PFO, we say that an edge $uv$ is fractionally oriented if $0<p(u,v)<1$. In that case we also say that the edge is \emph{unsaturated}.\footnote{This terminology is taken from a close relationship with network flow techniques.}  A cycle $C = u_1u_2\dots u_ku_1$ is said to be \emph{unsaturated} if all its edges $u_iu_{i+1}$ ($1\le i\le k$, where $u_{k+1}=u_1$) are unsaturated. For an unsaturated cycle we define $\alpha := \max\{p(u_i,u_{i+1})\mid 1\le i\le k\}$ and we note that $\beta:=1-\alpha>0$. By changing, for each $i\in[k]$, the value $p(u_i,u_{i+1})$ to $p(u_i,u_{i+1})+\beta$ and changing $p(u_{i+1},u_i)$ to $p(u_{i+1},u_i)-\beta$, we obtain another PFO in which all edges on the cycle, for which $p(u_i,u_{i+1})=\alpha$, become saturated, and so the cycle is no longer unsaturated. We will refer to this operation as \emph{saturating the cycle $C$}. By doing this, the corresponding potential outdegrees $d_p(u)$ of all vertices of the graph remain unchanged. It is easy to see that repeated use of this operation leads to a proof of Lemma~\ref{lem:rounding}.

We will need another operation that is similar to the one described above. Here we take a path $P = u_1\dots u_k$ ($k\ge2$), all of whose edges are unsaturated, and in addition to that, the endvertices $u_1$ and $u_k$ are also \emph{unsaturated}, meaning that none of $d_p(u_1)$ and $d_p(u_k)$ is an integer. Now we define $\alpha$ as the maximum of all values $p(u_i,u_{i+1})$ ($1\le i < k$) and set $\alpha' = d_p(u_1) - \lfloor d_p(u_1) \rfloor$ and $\beta' = d_p(u_k) - \lfloor d_p(u_k) \rfloor$. Similarly as above, we set $\beta = \min\{1-\alpha,1-\alpha',\beta'\}$ and then change each $p(u_i,u_{i+1})$ to $p(u_i,u_{i+1})+\beta$ and change each $p(u_{i+1},u_i)$ to $p(u_{i+1},u_i)-\beta$ ($1\le i < k$). This change produces new PFO in which either one of the edges on the path becomes saturated or one of the endvertices becomes saturated. We say that we have \emph{saturated the path $P$}. For further reference we state the basic property of the described change as the following fact.

\begin{observation}\label{obs:make path happy}
The process of saturating an unsaturated path $P=u_1\dots u_k$ described above produces new PFO $q$, in which either one of the edges on the path becomes saturated or one of the endvertices becomes saturated.
The potential degrees of all vertices except $u_1$ and $u_k$ remain unchanged, while $d_p(u_1)$ increases by $\beta$ and $d_p(u_k)$ decreases by $\beta$. Moreover, $d_q(u_1)\le \lceil d_p(u_1) \rceil$ and $d_q(u_k)\ge \lfloor d_p(u_k) \rfloor$.
\end{observation}

The well-known Hall Theorem characterizes when a bipartite graph with bipartition $U\cup V$ contains a matching covering all vertices in $U$. There is a weighted version in which a matching is replaced by an edge-set (subgraph) $M$ so that each vertex in $U$ is contained in the prescribed number of edges in $M$. We denote by $d_M(v)$ the number of edges in $M$ that are incident with the vertex $v$.
For $S\subseteq V(G)$, we write $N(S) = \cup_{s\in S} N(s)$, and for $v\in V$, we denote by $e(v,S)$ the number of edges from $v$ to $S$.
The following result by Ore \cite{Ore56,Ore57} is a generalization of Hall's Theorem, and is also called the Ore-Ryser Theorem by some authors (see \cite{VW13}).

\begin{lem}[Ore \cite{Ore56,Ore57}]
\label{lem:Generalized Hall}
Let $G=(U\cup V,E)$ be a bipartite graph and suppose that\/ $W: V\cup U \rightarrow \mathbb{N}_0$
is a weight function of the vertices of $G$. 
Suppose that for every vertex-set $S\subseteq U$ we have 
\begin{equation}
    \sum_{u\in S} W(u) \le \sum_{v\in V} \min \{W(v), e(v, S)\}.\label{eq:Hall condition}
\end{equation}
Then there is a subgraph $M$ of $G$ such that 
\begin{align*}
    d_M(u) &= W(u), \textrm{~~for each } u\in U, \textrm{ and } \\
    d_M(v) &\le W(v), \textrm{~~for each } v\in V.
\end{align*}
\end{lem}


\section{Main tools}
\label{sect:3}

In this section, we will introduce our main tools. Using them, we will be able to construct proper orientations in $r$-partite graphs. 

Throughout this section we will assume that $G$ is an $r$-partite graph ($r\ge2$) with the corresponding vertex-partition $V=V_1\cup \cdots \cup V_r$. We will also assume that $\MAD(G)\le 2k$, where $k$ is a positive integer, and we fix a $k$-orientation $D_0$ of $G$, and refer to $D_0$ as the \emph{base orientation} of $G$.

Suppose that $p$ is a partial fractional orientation (PFO) of $G$. Let us recall that the potential outdegree $d_p(v)$ of a vertex $v$ is defined as $d_p(v) = \sum_{u\in N(v)} p(v,u)$.
Note that $d_p(v)$ counts the number of unoriented edges incident with $v$ together with the outdegree $d_p^+(v)$ of the oriented edges. If all undirected edges incident with $v$ were oriented out of $v$, then (and only then) the outdegree would become equal to $d_p(v)$. We say that a vertex $v$ is \emph{oriented} if all its incident edges are (partially) oriented. 
On the other hand, if at least one edge incident with a vertex $v$ is unoriented, then we say that $v$ is an \emph{unoriented vertex}.

Equivalently, $d_p(v)=d_p^+(v)$.

We say that a PFO $p$ is \emph{$j$-proper} if the following holds:
\begin{enumerate}[parsep=0pt]
  \item[(1)] For every oriented vertex $v$, $d_p(v)$ is an integer.
  \item[(2)] For each $m>j$, the set $A_m$ of oriented vertices $v$ with $d_p(v)=m$ is independent.\footnote{A set of vertices is \emph{independent} if no two vertices in the set are adjacent.}
  \item[(3)] Every (partially) oriented edge is incident with a vertex in some $A_m$, where $m>j$.
  \item[(4)] If $u\notin \cup_{m>j}A_m$ and $v\in \cup_{m>j}A_m$ and $vu\in E$ is oriented out of $v$ in the base orientation, then $p(v,u)=1$ and $p(u,v)=0$, i.e. $vu$ is also oriented out of $v$ under $p$.

\end{enumerate} 
When referring to property (4) of a $j$-proper PFO, we will say that $p$ is \emph{aligned with $D_0$}.

Suppose that $p$ also satisfies:
\begin{enumerate}
  \item[(5)] (Strong version) If $v\in V(G)$ and $d_p(v)\le j$, then $d_p(v)$ is an integer.
\end{enumerate} 
Then $p$ is said to be \emph{strongly $j$-proper}.

We say that a $j$-proper PFO $p$ \emph{can be changed} into an $\ell$-proper PFO $q$, where $\ell\le j$, if the sets $A_m$ ($m>j$) are the same in $p$ and in $q$. Note that this implies (by (3)) that every oriented edge in $p$ is also oriented in $q$ (but $q$ may have more oriented edges if $\ell\ne j$).

By using the saturation process described at the end of the previous section, we are able to change any $j$-proper PFO into a strongly $j$-proper PFO.

\begin{lem}\label{lem:make strongly j-proper}
Suppose that a PFO $p$ is $j$-proper. Then $p$ can be changed into a strongly $j$-proper orientation $q$ such that the following holds:
\begin{itemize}
\item The set of unoriented edges is the same in $p$ and in $q$;
\item $d_q(v)=d_p(v)$ if $d_p(v)$ is an integer;
\item $d_q(v)\le d_p(v)$ if $d_p(v)>j$;
\item $d_q(v)$ is either equal to $\lfloor d_p(v)\rfloor$ or to $\lceil d_p(v)\rceil$ if $d_p(v)<j$.
\end{itemize}
\end{lem}

\begin{proof}
We change $p$ by making successive steps, in each step either decreasing the number of unsaturated edges with $0<p(u,v)<1$ or decreasing the number of vertices $v$ for which $d_p(v)<j$ is not integral. It suffices to describe one such step.

If for each vertex $v$ with $d_p(v)<j$, the value $d_p(v)$ is an integer, then $p$ is strongly $j$-proper and we stop. Otherwise, 
let $u$ be a vertex with noninteger value $d_p(u)<j$. Then $u$ is incident with an edge $uu_1$ with $0<p(u,u_1)<1$. If $d_p(u_1)$ is an integer, we find another edge $u_1u_2$ with $0<p(u_1,u_2)<1$ and continue building a path $uu_1u_2\dots u_t$ ($t\ge1$) until we either repeat one of the previous vertices ($u_t=u_i$ for some $1\le i< t$), or we reach a vertex $u_t$ with $d_p(u_t)$ non-integral, possibly $u_t=u$. In the latter case we have a path or a cycle $uu_1\dots u_t$ such that 
such that $d_p(u_i)$ is integer for $1\le i<t$ and $d_p(u_t)$ is not an integer. This path or cycle is clearly unsaturated. In the former case, we get an unsaturated cycle $u_iu_{i+1}\dots u_t$ such that $d_p(u_s)$ is integer for $i\le s<t$.

If we have an unsaturated cycle, then 
by saturating the cycle we keep potential degrees unchanged and we saturate at least one of the edges. If we have an unsaturated path, then by saturating this path, $d_p(u_t)$ decreases (but not more than down to the first smaller integer), so we are sure to have the third and the fourth property of the lemma, and $d_p(v)$ increases, but not more than to the ceiling of its value (see Observation \ref{obs:make path happy}), so we have the last property of the lemma. 

By doing this, either one of the edges on the path/cycle obtains value 0 or 1, or $d_p(v)$ or $d_p(u_t)$ becomes an integer. This completes the proof.
\end{proof}

Having the above lemma, we will be able to assume that every $j$-proper PFO is strongly $j$-proper.

With every (strongly) $j$-proper PFO $p$, we associate the \emph{gap values}. Let $A = \cup_{m>j} A_m$ be the set of oriented vertices in sets $A_m$ ($m>j$) from item (2) in the definition of $j$-proper. For a vertex $v\notin A$ we define its \emph{gap} with respect to $j$: 
$$gap(v) := d_p(v)-j.$$
If $G$ is $r$-partite with parts $V_1,\dots,V_r$, then we also define 
$$
     Gap(i) := \max \{ gap(v) \mid v\in V_i\setminus A \}, \quad 1\le i\le r.
$$
When speaking about gaps, the value $j$ will usually be clear from the context. If not, we may use the term \emph{$j$-gap} and use the notation $gap_j(v)$ and $Gap_j(i)$ to clarify.

\begin{lem}\label{lem:gap0}
Suppose that $G$ has a $j$-proper PFO $p$ such that each vertex that is not in $A = \cup_{m>j} A_m$ has nonpositive gap. Then $p$ can be changed into a proper orientation $q$ such that $d_q(v)\le \lceil d_p(v)\rceil$ for every $v\in V(G)$.

\end{lem}

\begin{proof}
By Lemma \ref{lem:make strongly j-proper}, we may assume that $d_p(v)$ is integer for each vertex $v$. 
Observe that every vertex with $d_p(v)>j$ is oriented since all gaps are nonpositive. Also, a vertex $u$ has gap zero precisely when $d_p(u)=j$.
Starting with an empty set $A_j$, add to $A_j$ any vertex $u$ with $gap_j(u)=0$ and orient all unoriented edges $uw$ incident with $u$ out of $u$, i.e. set $p(u,w)=1$ and $p(w,u)=0$. This will decrease the gap of $w$ by 1. 

Thus repeating this step will form an independent set $A_j$. After no vertex with gap 0 remains, all gaps are negative integers. Now, passing from $j$ to $j-1$, we have a $(j-1)$-proper PFO for which all gaps are nonpositive. Thus, we can repeat the process until all vertices have been oriented.

We end up with a 0-proper PFO. It remains to show that we can change all nonintegral $p$-values into 0/1, without changing the values $d_p(v)$. Let $H$ be the subgraph of $G$ containing those edges $uv$ for which $0<p(u,v)<1$. If there are no such edges, then we already have a proper orientation. Otherwise, $H$ contains a cycle $C$ since no vertex can be incident with precisely one edge in $H$. By saturating this cycle, we saturate one of the edges $uv\in E(C)$ and thus decrease the number of edges in $H$. By repeating this process, we obtain an orientation as claimed.
\end{proof}

In later sections we will construct proper orientations, starting with some $j$-proper partial orientation, and then we will decrease the gaps until no vertices with positive gap are left. We proceed from a $j$-proper PFO to a $(j-1)$-proper PFO such that the gaps $Gap(i)$ ($1\le i\le r$) decrease and eventually all become 0. Of course, in the new $(j-1)$-proper PFO, the gap is considered with respect to $j-1$. By doing so, the sets $A_m$ for $m>j$ will remain unchanged and we will form the new independent set $A_j$
If we achieve that $Gap(i)\le 0$ (with respect to the PFO being $(j-1)$-proper), we say that we \emph{close the gap} of the $i$th part $V_i$. And when all gaps become nonpositive, we say that the \emph{gaps were closed}. 

In our construction of proper orientations, we use the quantity $d_1(v)$, which denotes the number of edges which are in-edges for $v$ in the base orientation and are still unoriented under the partial fractional orientation we are building. The following simple corollary of properties (3) and (4) of $j$-proper PFOs will be used throughout.

\begin{lem}\label{lem:d1 bound}
Let $p$ be a $j$-proper PFO and $A = \cup_{m>j} A_m$. Then each vertex $v\notin A$ satisfies $d_p^+(v) \le k$ and
\begin{equation}
    d_1(v) \ge \ceilgapv + j - k.\label{eq:d1 bound}
\end{equation}
\end{lem}

\begin{proof}
Note that properties (3) and (4) of $j$-proper PFOs imply that each unoriented vertex $v$ has $d_p^+(v) \le k$ since each oriented edge $vu$ incident with $v$ is coming from some $A_m$, and if $p(v,u)>0$, then $vu$ is oriented from $v$ to $u$ in the base $k$-orientation $D_0$.  

Let $t$ be the number of unoriented edges incident with $v$ and let $t^+$ be the number of those edges that are outedges in $D_0$. Then  $d_p(v)=t+d_p^+(v) = (t-t^+) + (t^+ + d_p^+(v)) \le d_1(v) + k$.
Note that $gap(v) = d_p(v)-j$. 
Consequently, $d_1(v) \ge d_p(v) - k = gap(v) + j - k$.
Now, (\ref{eq:d1 bound}) follows since $d_1(v)$ is an integer.
\end{proof}

Let $l\ge r+k$ be an integer (the desired bound on the maximum outdegree of a proper orientation). We will construct an $(l-r)$-proper PFO $p$ such that $d_p(v)\le l$ for each $v\in V(G)$ by using the following algorithm. 

\bigskip
\noindent
{\bf Gap-Capping Algorithm}:\\[1mm]
We start with a partial orientation $p$ with all edges undirected. Now, we repeat the following process for $i=1,\dots,r$, in each step using the partial orientation $p$ obtained in previous steps:

\begin{itemize}
    \item[(i)] Let $X_{l-i+1}$ be the set of unoriented vertices with potential $d_p(v)\ge l-i+1$.
    \item[(ii)] Choose an independent set $A_{l-i+1}$ in $X_{l-i+1}$ containing all vertices in $X_{l-i+1}\cap V_i$ that has the maximum number of vertices subject to the condition that $X_{l-i+1}\cap V_i \subseteq A_{l-i+1}$.
    \item[(iii)] Orient each vertex $v\in A_{l-i+1}$ so that each $p$-unoriented edge $vu$ incident with $v$ that is oriented out of $v$ in $D_0$ is also oriented out from $v$ under $p$. Since $D_0$ is a $k$-orientation, this rule orients at most $k$ edges out of $v$. Since we had $d_p(v)\ge l-i+1 > k$ before this step, we can orient additional unoriented edges incident with $v$ such that $v$ becomes oriented and $d_p(v)$ becomes equal to $l-i+1$.
\end{itemize}

After each step of the algorithm, the partial orientation $p$ changes, and at the end it has properties which we summarize in the following lemma.

\begin{lem}\label{tool:rlevel}
If\/ $l\ge r+k$, then after performing the Gap-Capping Algorithm,
$p$ is a partial orientation of the graph $G$ and has the following properties:
\begin{itemize}
  \item[\rm (a)] $p$ is $(l-r)$-proper and is in particular aligned with $D_0$.
  \item[\rm (b)] $Gap_{l-r}(i)\le r-i$ for $1\le i\le r$.
\end{itemize}
\end{lem}

\begin{proof}
Let us first check properties (1)--(4) needed for an $(l-r)$-proper partial orientation. Property (1) holds since $p$ has no fractionally oriented edges. Similarly, (2) and (3) hold by the way we construct $p$ (we only orient edges incident to the vertices in sets $A_{l-i+1}$).

It remains to show that $p$ is aligned with $D_0$. Let $A = \cup_{m>l-r} A_m$. 
For a vertex $v\notin A$, suppose that $vu$ is an edge incident with $v$ that is oriented out of $v$. Then this edge was oriented when we have oriented the edges incident with $u$ and $u\in A_{l-i+1}$ for some $i$ ($1\le i\le r$). By (iii), the edge $vu$ is oriented from $v$ to $u$ in the base orientation $D_0$. This shows that $p$ is aligned with $D_0$.

It remains to prove (b). Consider an unoriented vertex $v\in V_i\setminus A$. When we defined $A_{l-i+1}$, $v$ was not included, so $d_p(v)$ was less than $l-i+1$ at that time. Since $d_p$ only decreases during the orientation process, we have the same condition also when the process is finished. This shows that $gap_{l-r}(v) = d_p(v) - (l-r) \le r-i$ and completes the proof.
\end{proof}

Starting with the partial orientation satisfying the statements of the above lemma, we will use our next lemma (the \emph{gap-decreasing tool}) which is based on the notion of a maximum-weight independent set and uses the Generalized Hall Theorem. In the proof of the lemma, we also make use of partial fractional orientations.

We will also use the following notation. For a vertex $v$ and $0\le i<r$, we define:
$$
 \delta_i(v) = \left\{
  \begin{array}{ll}
  0\,, & \hbox{if $gap(v)\le0$;} \\
  \frac{gap(v)}{d_1(v)-i\ceilgapv}, & \hbox{otherwise.}
  \end{array}
  \right.
$$
For a vertex-set $A$, we also write $\delta_i(A) = \max_{v\in A} \delta_i(v)$.

\begin{lem}\label{lem:tool3}
Let $p$ be a strongly $j$-proper PFO of an $r$-partite graph $G=(V_1,\dots,V_r,E)$, where $j\ge k$. Suppose that for some $i$, $0\le i< r$, at least $i$ gaps are nonpositive, say $Gap(2)\le0$, \ldots, $Gap(i+1)\le 0$. Suppose that for each unoriented vertex $v\in V_1$ with $gap(v)\ge 0$, $d_1(v)$ satisfies the following inequality:
\begin{equation}
d_1(v) \ge (i+1)\ceilgapv. \label{eq:condition1}
\end{equation}
Then we can change $p$ into a $(j-1)$-proper PFO and, meanwhile, for the $(j-1)$-gap, we close the gap of $V_1$ and keep the gaps of\/ $V_2,\dots,V_{i+1}$ nonpositive. For each other part $V_s$ $(s\ge i+2)$, if $Gap_{j-1}(s) > 0$, then 
$$
  Gap_{j-1}(s) \le Gap_j(s) + \max\left\{\delta_i(A_j\cap V_1), \delta_0(A_j\cap (V_{i+2}\cup\cdots \cup V_r)) \right\},
$$ 
where $A_j$ is the new independent set containing vertices with outdegree $j$, and the gap values in the definition of $\delta_i$ and $\delta_0$ refer to the $j$-gaps with respect to the original $j$-proper PFO $p$.
\end{lem}

\begin{proof}
In the proof we will define the color class $A_j$ by orienting some of the unoriented edges. First of all, let $A'$ be the set of oriented vertices whose gap is 0. For these vertices, property (3) of $j$-proper implies that all their neighbors are in $\cup_{m>j}A_m$. Thus, these vertices can be put in $A_j$ without worrying that they will be adjacent to any vertex in $A_j$. 
If $v$ is an oriented vertex with negative gap, then $gap_j(v)\le -1$ since $p$ is strongly $j$-proper. Thus, such a vertex will have $gap_{j-1}(v)\le0$, and thus we can henceforth neglect any vertices with negative gap.

Throughout the proof we will use the values $d_1(v)$. Recall that $d_1(v)$ is defined as the number of edges which are in-edges for $v$ in the base orientation and are still unoriented under the partial fractional orientation we are building. 
After orienting some of the edges, we change this value, and we will denote the value corresponding to the current orientation by $d_1'(v)$. Note that $d_1'(v) \le d_1(v)$.

Let $U$ be the set of all unoriented vertices with $gap(v)\ge 0$.
We define a weight function $W$ on $U$:
$$ W(v)=
    \begin{cases}
   \ceilgapv& v\in V_1, \\
    1 &   v\in V_s\ (2\le s\le i+1),\\
    0 &   otherwise.
    \end{cases}
$$
Next, we choose an independent set $A\subseteq U$ with maximum weight and, subject to this maximality condition, with $|A\cap V_1|$ as large as possible, and subject to these conditions, we also ask $|A|$ to be as large as possible. We let $X = U\setminus A$. Now, we claim that all vertices in $A$ can be oriented with outdegree $j$. From the definition of $W$, we obtain the following inequality for any set $S\subseteq X\cap V_s$ ($s\in\{1,\dots,i+1\}$):
\begin{equation}
    W(S)\le W(N(S)\cap A).\label{eq:Hall condition2}
\end{equation}

For each $s=2,\dots,i+1$, we now consider the bipartite graph $B_s$ with parts $X\cap V_s$ and $A\setminus V_s$ and all edges of $G$ between them. Having (\ref{eq:Hall condition2}), we can apply the Ore Theorem (Lemma \ref{lem:Generalized Hall}) to obtain an edge-set $M_s\subseteq E(B_s)$ such that 
\begin{align}
d_{M_s}(v) &\le W(v) \ \textrm{ if } v \in A\setminus V_s, \textrm{ and}\label{eq:Hall A-Vs}\\
d_{M_s}(v) &= W(v) \ \textrm{ if } v\in X\cap V_s.\label{eq:Hall X cap Vs}
\end{align}
Note that by (\ref{eq:Hall A-Vs}), $M_s$ does not have any vertex in $V_t$ for $t\ge i+2$.

Now, we orient the edges in $M_s$ from $A \setminus V_s$ to $X\cap V_s$ for $s=2,\dots,i+1$. After orienting these edges, each vertex in $X\cap (\bigcup_{s=2}^{i+1}V_s)$ gets at least one in-edge because of (\ref{eq:Hall X cap Vs}). So for every $v\in X\cap (\bigcup_{s=2}^{i+1}V_s)$, its gap decreases by at least 1. 
Next, since $Gap(s) \le 0$ (for $2\le s\le i+1$), we can orient all remaining unoriented edges incident with vertices in $A\cap(\bigcup_{s=2}^{i+1}V_s)$ out of this set. Since the gap of these vertices was 0, they will all end up being oriented with outdegree $j$. Moreover, each unoriented neighbor of any of these vertices will have its gap decreased at least by 1, since it gets at least one in-edge from $A\cap(\bigcup_{s=2}^{i+1}V_s)$. 

In particular, if $v\in X\cap V_1$, then $|N(v)\cap A\cap( V_{2}\cup \dots \cup V_{i+1})|\ge \ceilgapv+1$. Otherwise, we would have added the vertex $v$ into $A$ (and remove from $A$ the vertices in $N(v)\cap A\cap( V_{2}\cup \dots \cup V_{i+1})$), which would give the set of larger weight. So the $j$-gap of $v$ becomes smaller or equal to $-1$ (and therefore the $(j-1)$-gap will not be positive). 

Next, we describe how we orient the so far unoriented edges $vu$ incident with the vertices $v\in A\cap (V_1\cup V_{i+2}\cup\dots\cup V_r)$. First of all, we orient the out-edges $vu$ in the base orientation out of $v$ (so that their orientation coincides with that in $D_0$).
Note that this does not change the value of $d_1'(v)$.
Second, we use fractional orientation to let each $v\in A\cap (V_1\cup V_{i+2}\cup\dots\cup V_r)$ have outdegree $j$. For such a vertex $v$ and any unoriented edge $vu\in E$, we set 
\begin{equation}
    p(v,u) = 1-\frac{gap(v)}{d_1'(v)}\quad \textrm{and} \quad p(u,v) = 1-p(v,u).
    \label{eq:define p(v,u)}
\end{equation}

If $v\in A\cap V_1$, the assumption (\ref{eq:condition1}) and the orientation given in front, yield the following inequality:
$$d_1'(v)\ge (i+1)\lceil gap(v) \rceil-\sum_{2\le s\le i+1}d_{M_s}(v) \ge \lceil gap(v) \rceil.$$
The same holds for $v \in V_s\cap A$ $(s \ge i+2)$: we have $d_1'(v) = d_1(v) \ge \lceil gap(v) \rceil$ by (\ref{eq:d1 bound}) and the assumption that $j\ge k$.
This implies that $0 \le p(v,u)\le 1$ and thus the changed fractional orientation is well-defined. 
Each vertex $v\in A$ is now oriented and we claim that its outdegree is equal to $j$. We have already shown this for vertices in $A\cap(V_2\cup \dots \cup V_{i+1})$. For $V_1$ and for $V_s$, $s \ge i+2$, we had precisely $d_1'(v)$ unoriented edges incident with $v$ when we made the change (\ref{eq:define p(v,u)}). Thus, the outdegree of $v$ was increased precisely by $d_1'(v) - gap(v)$. Adding also the out-orientation of edges that are out-edges in the base orientation, this implies that the outdegree of $v$ is now precisely $j$.

Under the new orientation, the vertices in $A$ become oriented and have outdegree $j$. Meanwhile, any other unoriented vertex $u$, which is in $X\cap V_s$ $(s\ge i+2)$, gets at least one (fractionally) oriented edge since $A$ is of maximum size. So its $j$-gap decreases at least by $\min_{v\in A\cap N(u)}\{p(v,u)\}$. This means that the $(j-1)$-gap increases with respect to the $j$-gap by at most 
\begin{equation}
  \max_{v\in A\cap N(u)}\{1-p(v,u)\} 
  \le \max_{v\in A} \Bigl\{\frac{gap(v)}{d_1'(v)}\Bigr\}
  \le \max\left\{\delta_i(A_j\cap V_1), \delta_0(A_j\cap (V_{i+2}\cup\cdots \cup V_r)) \right\}, 
  \label{eq:max gap increase}
\end{equation}
with the last inequality holding because $d_1'(v) \ge d_1(v)-i\ceilgapv$ for $v\in A\cap V_1$ and because $d_1'(v) = d_1(v)$ for $v\in A\cap (V_{i+2}\cup\cdots \cup V_r))$.

Finally, we set $A_j := A\cup A'$ (where $A'$ is the set of vertices defined at the beginning of the proof). Note that the maximum in (\ref{eq:max gap increase}) can be taken over all of $A_j$ since for $v\in A'$, we have $gap(v)=0$. 

In conclusion, if initially we had $Gap_j(s)\ge0$ for $s \ge i+2$, then $Gap_{j-1}(s)$ may have increased, but the increase is at most the amount stated in the lemma. Of course, if $Gap_j(s)<0$, then $U\cap V_s=\emptyset$. In that case, the gap may increase by 1, but it will not become positive. 
This completes the proof.
\end{proof}

We have the following simplified corollary.

\begin{coro}\label{coro:tool4}
Let $p$ be a strongly $j$-proper PFO of an $r$-partite graph $G=(V_1,\dots,V_r,E)$. Suppose that $Gap(1)\ge0$ and that $Gap(2)\le0$, \ldots, $Gap(i+1)\le 0$ and that any unoriented vertex $v$ with $gap(v)\ge 0$ satisfies the following inequality:
\begin{equation}
    j\ge i\ceilgapv+k. \label{eq:condition Corollary 3.6}
\end{equation}
Then we can change $p$ into a $(j-1)$-proper PFO with a new independent set $A_j$ containing newly oriented vertices with outdegree $j$. Meanwhile, for the $(j-1)$-gap, we close the gap of $V_1$ and keep the gaps of $V_2,\dots,V_{i+1}$ nonpositive. For each other part $V_s$ $(s \ge i+2)$ with $Gap(s)\ge0$, the gap may increase but the increase cannot be bigger than  
$\max_{v\in A_j}\bigl\{\frac{gap(v)}{j-k-(i-1)\ceilgapv}\bigr\}$.
\end{coro}

\begin{proof}
Since $Gap(1)\ge0$, there is a vertex in $V_1$ whose gap is nonnegative. From (\ref{eq:condition Corollary 3.6}) we conclude that $j\ge k$. 
By (\ref{eq:d1 bound}) and (\ref{eq:condition Corollary 3.6}), we have $d_1(v) \ge \ceilgapv + j-k \ge (i+1)\ceilgapv$ for each $v\in V_1$. This shows that (\ref{eq:condition1}) holds, and we can apply Lemma \ref{lem:tool3}. The lemma gives a $(j-1)$-proper PFO and the new independent set $A_j$. If $v\in A_j$, then (\ref{eq:d1 bound}) implies that $j-k-(i-1)\ceilgapv \le d_1 - i\ceilgapv$. Therefore, for $v\in A_j\cap V_1$ with $gap(v)>0$,
$$
   \frac{gap(v)}{j-k-(i-1)\ceilgapv} \ge \frac{gap(v)}{d_1(v)-i\ceilgapv}.
$$
Similarly, if $v\in A_j\cap V_s$ ($i+2\le s\le r$) has $gap(v)>0$, then by (\ref{eq:d1 bound}) and (\ref{eq:condition Corollary 3.6}), we have $d_1(v) \ge \ceilgapv + j-k \ge (i+1)\ceilgapv$. So, we have for every $i\ge0$: 
$$
   \frac{gap(v)}{j-k-(i-1)\ceilgapv} \ge \frac{gap(v)}{j-k+\ceilgapv} \ge \frac{gap(v)}{d_1(v)}.
$$
These inequalities combined with Lemma \ref{lem:tool3} confirm that $\max\{\delta_i(A_j\cap V_1), \delta_0(A_j$ $\cap  (V_{i+2}\cup\cdots \cup V_r)) \}\le\frac{gap(v)}{j-k-(i-1)\ceilgapv}$. So the last conclusion of the corollary holds.
\end{proof}

\section{Proper orientations of bipartite graphs}

In this section we prove our first main result, Theorem \ref{th:main3}, which gives a bound for the proper orientation number of bipartite graphs. The proof is split into two parts. First we establish the bounds, and then we show that the bounds are best possible.

First we prove that for every bipartite graph $G$, 
\begin{equation}
  \bigl\lceil\tfrac{1}{2}\MAD(G)\bigr\rceil \le \Vec{\chi}(G) \le \bigl\lceil\tfrac{1}{2}\MAD(G)\bigr\rceil+3.
   \label{eq:Upper bound to be proved}
\end{equation}

\begin{proof}[Proof of (\ref{eq:Upper bound to be proved}).]
The lower bound in (\ref{eq:Upper bound to be proved}) holds for any graph. To prove the upper bound, let $G=(V_1,V_2,E)$ be a bipartite graph with $\MAD(G)\le 2k$ and let $l=k+3$. As discussed in the previous section, we start with a base orientation $D_0$, which is a $k$-orientation of $G$. Then we use Lemma \ref{tool:rlevel} to obtain an $(l-2)$-proper PFO $p$ such that $Gap(1)\le 1$ and $Gap(2)\le 0$. Since $l-2\ge k+1 \ge \ceilgapv+k$ for $v\in V_1$, we can apply Corollary~\ref{coro:tool4} (by using $i=1$). By doing this, we determine $A_{l-2}$ and close all gaps. Thus we can change $p$ into an $(l-3)$-proper PFO with all gaps nonpositive. By Lemma \ref{lem:gap0}, there is a proper orientation $q$ for $G$ such that $d_q(v)\le d_p(v)$. So $d_q(v)=d_q^+(v)\le l = k+3$ for each vertex. This shows that $\Vec{\chi}(G)\le k+3$. 
\end{proof}

\subsubsection*{Tightness of the upper bound}

Next we prove that there exist a bipartite graph $G$ such that 
$\Vec{\chi}(G) = \bigl\lceil\tfrac{1}{2}\MAD(G)\bigr\rceil+3$.

We first define a bipartite graph $G_1$ with the bipartition $\{A\cup D, B_1\cup \cdots \cup B_k \cup C\}$. The set $A$ has $k$ vertices $v_1,\dots,v_k$.  
For $i\in[k]$, the set $B_i$ has $(k(k+2)+1)\binom{k}{i}$ vertices. For each $i$-subset $S\subseteq A$, a set of $k(k+2)+1$ of the vertices in $B_i$ is adjacent precisely to the vertices in $S$. The set $C$ has $m = (k+3)k^2$ vertices and is completely joined to $A$. Divide $C$ into pairwise disjoint $k$-subset. Finally, the set $D$ has $t=\tfrac{m}{k}$ vertices, each of which is joined to a different $k$-set of the partition of $C$. It is not hard to see that $G_1$ is $k$-degenerate, i.e., it can be reduced to the empty graph by successively deleting vertices of degree at most $k$. This implies that $\MAD(G_1)< 2k$. 
 
\begin{claim}\label{claim0}
Suppose that $p$ is an orientation of $G_1$ such that $d_p^+(v) \le k+2$ for every vertex $v$ and that $d_p^+(x)\ne d_p^+(y)$ for every edge $xy$, where $x\in D$ and $y\in C$. Then for each $i\in [k]$ and each $v\in A$, there is a vertex in $B_i$ with outdegree $i$ that is adjacent to $v$. Moreover, $C$ contains a vertex with outdegree $k+1$.
\end{claim}

\begin{proof}
  Take $S\subseteq A$ with $v\in S$ and $|S|=i$, where $i\in [k]$. If each vertex in $B_i$ that is completely adjacent to $S$ (i.e.\ adjacent to each vertex in $S$) has outdegree less than $i$, then the number of edges oriented from $S$ to $B_i$ is at least $k(k+2)+1$, which implies that some vertex in $S$ has outdegree more than $k+2$. This contradiction proves the first claim.
  
  Suppose now that no vertex in $C$ has outdegree $k+1$. Since $d_p^+(v)\le k+2$ for each $v\in A$, there are at most $k(k+2)$ vertices in $C$ that have an incoming edge from $A$. Since there are $(k+3)k$ pairwise disjoint $k$-subsets in $C$, one of the parts, say $X$, of the partition of $C$ into $k$-sets has all edges directed from $X$ to $A$. The vertices in $X$ thus have outdegree $k$ towards $A$. If all of them have outdegree in $G_1$ exactly $k$, their common neighbor in $D$ would have outdegree exactly $k$ as well, and this would not be a proper orientation.
  Since each vertex in $X$ has degree $k+1$, and outdegree more than $k$, its outdegree is precisely $k+1$.
\end{proof}

Now, we take 8 disjoint copies of $G_1$, denoted by $G_1,\dots,G_{8}$. We denote the vertices in the set $A$ of $G_s$ by $A^{(s)}$ $(1\le s\le 8)$. Finally, add all edges between $A^{(i)}$ and $A^{(j)}$ for all pairs $(i,j)\in \{(1,4),(2,4),(3,4),(4,5),(5,6),(5,7),(5,8)\}$ and denote the resulting graph by $G$.

\begin{claim}\label{claim1}
The graph $G$ defined above is bipartite and has 
$\MAD(G) < 2k$. Moreover, $\Vec{\chi}(G) = k+3$.
\end{claim}

\begin{proof}
It is not hard to see that $G$ is bipartite and that 
$\MAD(G) < 2k$ (since $G$ is $k$-degenerate).
By (\ref{eq:Upper bound to be proved}) it is also clear that $\Vec{\chi}(G) \le k+3$.

In order to show that $\Vec{\chi}(G) = k+3$, let us suppose, for a contradiction, that there is a proper orientation $p$ of $G$ such that $d_p^+(v)\le k+2$ for every vertex $v$. Claim \ref{claim0} implies that for each $s\in [8]$ and each $v\in A^{(s)}$, we have $d_p^+(v) = 0$ or $k+2$. Then it is easy to see that all vertices in $A^{(1)}\cup A^{(2)}\cup A^{(3)}\cup A^{(5)}$ have the same outdegree (either $0$ or $k+2$), and that all vertices in $A^{(4)}\cup A^{(6)}\cup A^{(7)}\cup A^{(8)}$ have the same outdegree (which is different from the outdegrees in $A^{(5)}$). Without loss of generality, we may assume that the first set has common outdegree 0. However, this implies that the vertices in $A^{(4)}$ have outdegree at least $4k$, which is larger than $k+2$, a contradiction.
\end{proof}

\section{Proper orientations of $r$-partite graphs}
\label{sect:r-partite}

In this section, we use Corollary~\ref{coro:tool4} to give a proof of Theorem \ref{th:main2}. In fact, we will prove the following more specific result.

\begin{theo}\label{thm:explicit upper bound}
Let\/ $G$ be an $r$-partite graph with $\MAD(G)\le 2k$. Let $t$ be the smallest integer for which $t^{t+1} \ge r-1$. Then
$$
  \Vec{\chi}(G) \le k + 3r(t+1) \le k + \frac{3(1+o(1))\, r\log r}{\log\log r}.
$$
\end{theo}

\begin{proof}
Let $\chi(G)=r$, $\MAD(G)\le 2k$ and $l=k+3r(t+1)$.
We will prove that there exists a proper $l$-orientation $p$ of $G$. First, we will deal with vertices whose degree is large. As proved before, this can be controlled by starting with a base $k$-orientation $D_0$ and use Lemma~\ref{tool:rlevel}. Having done that, we have obtained an $(l-r)$-proper PFO $p$ such that $Gap(s) \le r-s$ for each $s=1,\dots,r$. 

For convenience, let $V_{i}=V_{i \bmod r}$ for any positive integer $i$. We will use Corollary~\ref{coro:tool4} repeatedly and will determine $A_{l-r},\dots,A_{l - r(t+1) + 1}$ iteratively. 
In order to apply Corollary~\ref{coro:tool4} in this iteration and also in later steps, we have to make sure that the ``$j\ge k$'' condition from the corollary is satisfied. To see this, note that $l=k+3r(t+1)$, so in our process, $j\ge l-r(t+1)+1\ge k$. 

When we determine $A_{l-r-i}$ $(0\le i\le rt-1)$, we view $V_{i+1}$ as the set $V_1$ in Corollary~\ref{coro:tool4} (which is used throughout with its value of $i=0$, i.e., we do not insist on any parts keeping their nonpositive gap). When applying the corollary, we close the gap of $V_{i+1}$, and gaps of other parts increase at most by 
$\max_{v\in A_{l-r-i}}\{\frac{gap(v)}{(l-r-i)-k}\}$. 
For this process, we have following claim.

\begin{claim}
Let $j=\lceil \frac{i+1}{r}\rceil$. When we determine $A_{l-r-i}$ $(0\le i\le rt-1)$, we close the gap of $V_{i+1}$ and for any other part, if the gap is positive, then it has increased by at most $t^{-j}$ in this step.
\end{claim}

\begin{proof}
From Corollary \ref{coro:tool4}, it is sufficient to prove that $\frac{gap(v)}{(l-r-i)-k}\le t^{-j}$ for every $v\in A_{l-r-i}$.

The proof is by induction on $i$. If $i=0$, then for every vertex $v$, we have:
\begin{equation}
    \frac{gap(v)}{(l-r-i)-k}\le \frac{r-1}{l-r-k} \le \frac{r}{rt} = \frac{1}{t}.
    \label{eq:i=0}
\end{equation}
Through the next steps when $i=1,\dots,r-1$, the induction hypothesis shows that the gap of each vertex is at most $r-1$, and hence (\ref{eq:i=0}) holds. Thus, we may now assume that 
$i\ge r$ and $j\ge2$. Since the gap of each part was nonpositive at most $r-1$ steps earlier, the induction hypothesis implies that the gap of each vertex is currently at most $(r-1)t^{1-j}$. Thus, $$\frac{gap(v)}{(l-r-i)-k}\le \frac{(r-1)\,t^{1-j}}{(l-r-rt+1)-k} \le \frac{rt^{1-j}}{tr} = t^{-j}.$$
\end{proof}

After finishing the above process iteratively for $i=0,1,\dots,rt-1$, we have a strongly $(l-r(t+1))$-proper PFO $p$, whose gaps are bounded: $Gap(x)\le \frac{r-x}{t^t}$ for $1\le x\le r$. Next, we use Corollary \ref{coro:tool4} iteratively to determine $A_{l-r(t+1)},\dots, A_{l-r(t+2)+1}$ and close the gap of each part. The following claim will help us proving that this process can be made so that no positive gaps remain.

\begin{claim}
When we determine $A_{l-r(t+1)-i}$ $(i=0,\dots,r-1)$, we close the gap of $V_{i+1}$, and keep the gaps of\/ $V_1,\dots,V_i$ nonpositive. For any other part, the gap increases by at most $t^{-(t+1)}$.
\end{claim}

\begin{proof}
The proof is by induction on $i$. When $i=0$, before using Corollary~\ref{coro:tool4}, we need to check the inequality in Corollary~\ref{coro:tool4}. For $v\in V(G)$ with $gap(v)\ge 0$, $l-r(t+1)\ge k$.
Then for a vertex $v$ in $V_2, \dots, V_r$, we have $gap(v)\le (r-1)t^{-t}$ and hence
$$\frac {gap(v)}{l-r(t+1)-k} \le 
\frac{(r-1)t^{-t}}{l-r(t+1)-k} \le \frac{rt^{-t}}{rt} = t^{-(t+1)}.$$
This shows that the gaps of $V_{2}, \dots, V_r$ increase by at most $t^{-(t+1)}$ (or remain nonpositive). 
Next, when $i>0$, for a vertex in $V_{s}$ ($i+2\le s\le r$) we have $gap(v)\le \frac{r-i}{t^t}+\frac{i-1}{t^{(t+1)}}\le \frac{r-1}{t^t}\le t$ by the induction assumption. Then we check the inequality in Corollary~\ref{coro:tool4}:
$$(i-1)\ceilgapv+k\le tr+k \le l-r(t+1)-i.$$
And 
$$\frac {gap(v)}{l-r(t+1)-i-k-i\ceilgapv} \le 
\frac{rt^{-t}}{l-r(t+1)-r-k-rt} \le \frac{rt^{-t}}{rt} = \frac{1}{t^{(t+1)}}.$$
This shows that the gaps of $V_{i+2}, \dots, V_r$ increase by at most $\frac{1}{t^{(t+1)}}$ and completes the proof of the claim.
\end{proof}

After determining $A_{l-r(t+2)+1}$, we have $Gap(s)\le 0$ for $s=1,\dots,r$. By Lemma \ref{lem:gap0}, there is a proper orientation for $G$, whose maximum outdegree is at most $l$. 
\end{proof}

\section{Proper orientations of planar graphs}

Results of the previous section show that the proper orientation number of planar graphs is bounded. 
Theorem \ref{th:main2} combined with the fact that planar graphs are 4-colorable yields a rough upper bound. In this section, we will use our main tools together with a more detailed analysis to obtain a better bound. Our final result, Theorem \ref{th:main1}, will then establish that for every planar graph $G$, we have $\Vec{\chi}(G)\le 14$. Lastly, we will also outline the proof that every outerplanar graph $G$ has $\Vec{\chi}(G)\le 10$.

\begin{proof}[Proof of Theorem \ref{th:main1}]
Let $G$ be a planar graph. Let us first observe that $\MAD(G)<6$. This well-known fact follows easily from Euler's formula. Thus $G$ has a 3-orientation $D_0$, which we call the base orientation as before. Furthermore, by the Four-Color Theorem \cite{AH89,RSST97}, $G$ is 4-partite. We let $V_1,\dots,V_4$ be the color classes of a 4-coloring of $G$.

Next we give a construction of a proper orientation $p$. We let $l=14$ as an upper bound on the desired maximum outdegree. 

\medskip
\noindent
{\bf Step 0.}
By Lemma~\ref{tool:rlevel}, we construct an $(l-4)$-proper partial orientation $p$ such that $Gap(i)\le 4-i$ $(i=1,2,3,4)$. Moreover, $p$ is aligned with the base orientation $D_0$. This means that any vertex $v\notin A_l\cup A_{l-1}\cup A_{l-2}\cup A_{l-3}$ has $d_p^+(v)\le 3$. 

\medskip
Our next goal is to close all gaps. We achieve this by using Lemma~$\ref{lem:tool3}$ in four additional steps:

\medskip
\noindent
{\bf Step 1.}
We determine $A_{l-4}$, making $p$ an $(l-5)$-proper PFO, closing $Gap(1)$, keeping $Gap(4)$ nonpositive, and achieving $Gap(2)\le 2+3/7$ and $Gap(3)\le 1+3/7$. 

\medskip
\noindent
{\bf Step 2.}
We determine $A_{l-5}$, changing $p$ into an $(l-6)$-proper PFO, closing $Gap(2)$, keeping $Gap(1)$ and $Gap(4)$ nonpositive, and achieving $Gap(3)\le 1+3/7+17/24$. 

\medskip
\noindent
{\bf Step 3.}
We determine $A_{l-6}$, making $p$ an $(l-7)$-proper PFO, closing $Gap(3)$, keeping $Gap(1)$ and $Gap(2)$ nonpositive, and achieving $Gap(4)\le 1$. 

\medskip
\noindent
{\bf Step 4.}
Finally, we change $p$ into an $(l-8)$-proper PFO and close all gaps.

\medskip

Let us now describe how to achieve the claimed statements in Steps 1--4. Having the $(l-4)$-proper partial orientation from Step 0, we will find $A_{l-4}$, close the gap of $V_1$ and keep $Gap(4)$ nonpositive by applying Lemma \ref{lem:tool3} with $i=1$. We need to confirm the inequality (\ref{eq:condition1}). Let $v$ be an unoriented vertex. Since $p$ is aligned with $D_0$, we have:
\begin{equation}
  d_1(v) \ge d_p(v) - 3 = (gap(v)+l-4)-3 = gap(v)+7 \ge 2 \ceilgapv.\label{eq:d1 after step 0}
\end{equation}
Note that all gaps are integers at this stage.
By Lemma \ref{lem:tool3} we now obtain an $(l-5)$-proper PFO and close the gap for $V_1$. By using (\ref{eq:d1 after step 0}), we can conclude that the gaps of $V_2$ and $V_3$ increase by at most 
$$
  \max\left\{\delta_1(A_{l-4}\cap V_1), \delta_0(A_j\cap (V_{2}\cup V_{3})) \right\}\le \max_{v\in A_{l-4}}\Bigl\{\frac{gap(v)}{d_1(v)-\ceilgapv}\Bigr\}.
$$
This yields the following upper bound on the increase of the gaps: 
$$
  \max_{v\in A_{l-4}}\Bigl\{\frac{gap(v)}{d_1(v)-\ceilgapv}\Bigr\} \le
  \max_{v\in A_{l-4}}\Bigl\{\frac{gap(v)}{l-7}\Bigr\} \le \frac{3}{7}.
$$

We conclude that the new gaps satisfy $Gap(2)\le 2+\frac{3}{7}$ and $Gap(3)\le 1+\frac{3}{7}$. 

Let us also observe that any vertex $u\in V_2$ with its new value $gap(u)>2$ belongs to $V_2$, had its $(l-4)$-gap equal to 2 and in the process of changing $p$ to become $(l-5)$-proper, precisely one edge $vu$ incident with $u$ was (fractionally) oriented and the edge $uv$ is directed out of $u$ in the base orientation. This means that $d_1(u)$ remains the same. In particular, we have $d_1(u)\ge 9$ by (\ref{eq:d1 after step 0}).

In Step 2, we determine $A_{l-5}$, close the gap of $V_2$ and keep the gaps of $V_1$ and $V_4$ nonpositive. Here we take $i=2$ while applying Lemma \ref{lem:tool3} (with appropriate permutation of $V_1,V_2,V_3,V_4$). As remarked above, if $gap(u)>2$, then we have:
$$
   d_1(u)\ge 9 = 3\lceil gap(u)\rceil.
$$
If $gap(u)\le 2$, then 
$$d_1(u)\ge \lceil gap(u)\rceil + (l-5)-3 = \lceil gap(u)\rceil + 6 \ge 3\lceil gap(u)\rceil.$$

So we can apply Lemma \ref{lem:tool3} with $i=2$ to obtain $A_{l-5}$ and change $p$ into an $(l-6)$-proper PFO, while closing the gap of $V_2$ and keeping gaps of $V_1$ and $V_4$ nonpositive. For $V_3$, the gap increases at most by 
$$
  \max\left\{\delta_2(A_{l-5}\cap V_2), \delta_0(A_j\cap V_{3}) \right\} \le
\max_{v\in A_{l-5}} \Bigl\{\frac{gap(v)}{d_1(v)-2\ceilgapv} \Bigr\}.
$$
If $gap(v)>2$, then as shown above, $d_1(v) - 2\ceilgapv = d_1(v)-6 \ge 3$. Thus we have:
$$\frac{gap(v)}{d_1(v)-2\ceilgapv}\le \frac{2+\tfrac{3}{7}}{3} \le \frac{17}{21}.$$
If $gap(v)\le 2$, then $d_1(v) \ge 6 + gap(v)$ and we have:
$$
  \frac{gap(v)}{d_1(v)-2\ceilgapv}\le \frac{gap(v)}{6+gap(v)-2\ceilgapv}\le \frac{gap(v)}{2+gap(v)}\le \frac{1}{2}.
$$
We conclude that $Gap(3)\le 1+\frac{3}{7}+\frac{17}{21}$ and that $Gap(1),Gap(2),Gap(4)\le 0$. 

In Step 3 we will close the gap of $V_3$ and keep the gaps of $V_1$ and $V_2$ nonpositive. In order to achieve this, we take $i=2$ and apply Lemma~\ref{lem:tool3}. And we see if $gap(v)>2$, $v$ has at least two fractionally oriented edges in the current $(l-6)$-proper PFO~$p$. These edges are both out-edges for $v$ in $D_0$ and we see in the same way as above that $d_1(v) \ge 9 = 3 \ceilgapv$. Similarly, if $1 < gap(v) \le2$, $v$ gets at least one fractionally oriented edge in $p$ and this edge is an out-edge for $v$ in $D_0$. Consequently, 
$$d_1(v)\ge (\ceilgapv+l-5)-3 = 8 \ge 3\ceilgapv.$$ 
If $gap(v)\le 1$, then $d_1(v)\ge (\ceilgapv+l-6)-3 \ge 3\ceilgapv$.

We conclude that the condition (\ref{eq:condition1}) holds and thus we
obtain an $(l-7)$-proper PFO and close the gap for $V_3$. When applying Lemma \ref{lem:tool3}, the gap of $V_4$ increases at most by
$$\max_{v\in A_{l-6}} \Bigl\{\frac{gap(v)}{d_1(v)-2\ceilgapv}\Bigr\}.$$

Consider any vertex $v\in A_{l-6}$. If $gap(v)>2$, then
$$\frac{gap(v)}{d_1(v)-2\ceilgapv}\le \frac{gap(v)}{9-2\ceilgapv}\le \frac{gap(v)}{3}\le 1.$$
If $gap(v)\le 2$, then
$$\frac{gap(v)}{d_1(v)-2\ceilgapv}\le \frac{gap(v)}{4-\ceilgapv}\le \frac{gap(v)}{2}\le 1.$$
This implies that the new gaps are bounded as follows: $Gap(4)\le 1$ and $Gap(1), Gap(2),$ $ Gap(3) \le 0$. 

In the last step, we close the gap of $V_4$ and keep the gaps of $V_1,V_2,V_3$ nonnegative. We only need to check the inequalities to apply Lemma \ref{lem:tool3} (with $i=3$).
If $gap(v) > 0$, then $v\in V_4$ and $v$ gets at least one fractionally oriented edge in Step 3. So 
$$d_1(v)\ge (\ceilgapv +l-8)-2 \ge 4 = 4\ceilgapv.$$ 
And if $gap(v) = 0$, we have $d_1(v) \ge (l-8)-3\ge 0 = 4\ceilgapv$. 
Thus we can use Lemma~\ref{lem:tool3} to close the gap of $V_3$, while keeping the other gaps nonpositive. Finally, after no gaps are positive, Lemma~\ref{lem:gap0} shows that $p$ can be changed to a proper orientation of $G$. As we started with $l=14$, we get the claim of the theorem.
\end{proof}

\begin{proof}[Proof of Theorem \ref{th:planar3colorable}]
For 3-colorable graphs, we can use the same proof as above, except that we start with three classes $V_1,V_2,V_3$. Considering a graph $G$, we use its base $3$-orientation $D_0$, and if $G$ is outerplanar, we can assume that $D_0$ is a 2-orientation. In order to treat both cases at the same time we set $l=11$ and $k=3$ for the general 3-colorable case, and $l=10$ and $k=2$ for the outerplanar case. 

Next we sketch how to construct a proper orientation $p$ with maximum outdegree at most $l$. By Lemma~\ref{tool:rlevel}, we can construct an $(l-3)$-proper partial orientation $p$ such that $Gap(i)\le 3-i$ for $i=1,2,3$. Meanwhile, for any unoriented vertex $v\in V_i$, $d_p^+(v)\le k$, and the partial orientation is aligned with $D_0$. 
Then we need to close the gap of these parts. We achieve this by using Lemma~$\ref{lem:tool3}$ in two steps:

\medskip
\noindent
{\bf Step 1.}
We determine $A_{l-3}$, making $p$ $(l-4)$-proper, closing $Gap(1)$, keeping $Gap(3)$ nonpositive, and achieving $Gap(2)\le 1+\frac{1}{2}$. 

\medskip
\noindent
{\bf Step 2.}
We determine $A_{l-4}$, making $p$ $(l-5)$-proper, closing all gaps.

\medskip

Next, we describe how to achieve the claimed statements. First we find $A_{l-3}$, close the gap of $V_1$ and keep $Gap(3)$ being nonpositive. We check the inequlities in Lemma \ref{lem:tool3} where 
$$d_1(v)\ge (gap(v)+l-3)-k\ge gap(v)+5\ge 2\ceilgapv.$$
So we obtain an $(l-4)$-proper PFO and close the gap of $V_1$. And for $V_2$, the gap increases at most by 

$$
  \max_{v\in A_{l-3}} \Bigl\{\frac{gap(v)}{d_1(v)-\ceilgapv}\Bigr\} \le 
  \max_{v\in A_{l-3}} \Bigl\{\frac{gap(v)}{gap(v)+5-\ceilgapv}\Bigr\}
  \le \frac{2}{5}.
$$ 
Here we have used the fact that $gap(v)\le2$ for every $v\in A_{l-3}$.

In Step 2 we have $Gap(2)\le 1+\frac{2}{5}$ and $Gap(1),Gap(3)\le 0$. Next we find $A_{l-4}$ and close all gaps. To check the inequality needed in order to apply Lemma \ref{lem:tool3}, we first consider any vertex $v$ with $gap(v) > 1$. Clearly, $v\in V_2$, and $v$ gets at least one fractionally oriented edge in Step 1. Such an edge is an out-edge for $v$ in the base orientation. Thus, if $1<gap(v)<2$, then 
$$d_1(v)\ge (\ceilgapv +l-4) -2 \ge 6 = 3\ceilgapv.$$
If $gap(v)\le 1$, then 
$$d_1(v)\ge gap(v)+l-4- 2\ge 4+gap(v)\ge  3\ceilgapv.$$
Thus we obtain an $(l-5)$-proper PFO and $Gap(1),Gap(2),Gap(3)\le 0$. By Lemma~\ref{lem:gap0}, we can now obtain an $l$-proper orientation. So $\Vec{\chi}(G) \le l=11$.
\end{proof}

\section{Conclusion}

In this paper, we introduced a tool combining fractional orientations and Hall theorem. By repeated use of it, we have proved that $\Vec{\chi}(G)$ can be bounded by a function of $\MAD(G)$ and $\chi(G)$. In this way, we also gave detailed upper bounds for planar and outerplanar graphs, by using the fact that every planar graph is 4-colorable with $\MAD(G)$ at most 3, and that every outerplanar graph is 3-colorable with $\MAD(G)$ at most 2. But, there are still some basic questions about proper orientation number that remain open.

The first one is about the tightness of our bound on proper orientation number of planar graphs.

\begin{ques}
What is the maximum value of $\Vec{\chi}(G)$ over all planar graphs $G$?
\end{ques}

We know that the answer to this question is between 10 and 14, but it is unclear what is the right number.

In Corollary \ref{cor:main2}, we show that $\Vec{\chi}(G)$ is bounded by a (non-linear) function of $\MAD(G)$. It would be important to know whether this function can be replaced by a linear upper bound.

\begin{ques}
Can $\Vec{\chi}(G)$ be bounded above by a linear function of\/ $\MAD(G)$?
\end{ques}

Theorem \ref{th:main2} shows that $\Vec{\chi}(G)$ is bounded by the sum of $\tfrac{1}{2}\MAD(G)$ and a non-linear function of $\chi(G)$. As $\frac 12 \MAD(G)$ is a trivial lower bound on $\Vec{\chi}(G)$, the gap between them is more interesting. 

\begin{ques}
Can $\Vec{\chi}(G)-\frac12 \MAD(G)$ be bounded by a linear function of $\chi(G)$?
\end{ques}

In particular, for a graph $G$ of genus $g$, Corollary \ref{cor:main2} implies that 
$$\Vec{\chi}(G) \le O(g^{1/2}\log g / \log\log g).$$
Is there a better bound? We believe there is and propose the following.

\begin{conj}
If $G$ is a graph of genus $g$, then $\Vec{\chi}(G) \le O(g^{1/2})$.
\end{conj}

\bibliographystyle{abbrv}
\bibliography{biblio_proper_orientations}

\end{document}